\documentclass{article}
\author{Romain Bondil}
\title{General elements of an $m$-primary ideal on a normal surface singularity}

\usepackage{amsmath,amssymb,amsthm,amsfonts,latexsym}

\newtheorem{defi}{Definition}[section]
\newtheorem{thm}[defi]{Theorem}
\newtheorem{lem}[defi]{Lemma}
\newtheorem{prop}[defi]{Proposition}
\newtheorem{cor}[defi]{Corollary}
\theoremstyle{definition}
\newtheorem{rem}[defi]{Remark}

\newcommand{\C}{\mathbb{C}}

\newcommand{\Dc}{\mathcal{D}}

\newcommand{\Oc}{\mathcal{O}}

\newcommand{\OSO}{\Oc_{S,0}}

\newcommand{\SIB}{\overline{S_I}}

\newcommand{\XIB}{\overline{X_I}}
\newcommand{\OXIB}{\mathcal{O}_{\XIB}}
\newcommand{\bIb}{\overline{b_I}}

\DeclareMathOperator{\dgre}{deg}

\def\DebEq{\vskip-8pt}
\def\FinEq{\vskip-3pt \noindent}

\begin{document}
\title{General elements of an $m$-primary ideal on a normal surface singularity}
\author{Romain Bondil}

\date{}

\maketitle

\begin{abstract}
In this paper\footnote{M.S.C. 32S15, 32S25 and 14J17, 14H20.}, we show how to apply a theorem by L\^e D.T. and the author about linear families of curves on normal surface singularities  to get new results in this area. The main concept used is a specific definition of {\em general elements} of an ideal in the local ring of the surface. We make explicit the connection between this notion and the elementary notion of general element of a linear pencil, through the use of {\em reduction}. This allows us to prove the invariance of the generic Milnor number (resp. of the multiplicity of the discriminant), between two pencils generating two ideals with the same integral closure (resp. the  projections associated). We also show that our theorem, applied in two special cases, on the one hand completes a theorem by Snoussi on the limits of tangent hyperplanes, and on the other hand gives an algebraic $\mu$-constant theorem in linear families of  planes curves.
\end{abstract}


%

\section*{Introduction}

Let $(S,0)$ be a germ of normal complex-analytic surface, with local ring $\OSO$ corresponding to the germs of holomorphic functions on $(S,0)$, and maximal ideal $m$, formed by the germs taking  the value $0$ at $0$.

To  any couple $(f,g)$ of elements of $m$, one may associate three related objects~: the  {\em linear pencil} of the curves $C_{\alpha,\beta} : \alpha f + \beta g=0$ with $(\alpha,\beta)\in \C^2$, the  {\em ideal}  $J=(f,g)$ in $\OSO$, and the  {\em projection}~: 
\DebEq \begin{eqnarray*}
p : (S,0)\rightarrow (\C^2,0), \\
x \mapsto (f(x),g(x)).
\end{eqnarray*}
\FinEq
We will always assume that the curves $f=0$ and $g=0$ share no common component (in other words~: the corresponding linear system has no  {\em fixed component}, the ideal $J$ is {\em $m$-primary}, and   the projection $p$ is {\em finite}).

Denoting  by $(\Delta_p,0)\subset (\C^2,0)$ the discriminant of the projection $p$ (see \S~\ref{sec-discri}), one may define a {\em general element of the pencil} $(C_{\alpha,\beta})$ as  the inverse-image by~$p$ of any  line $\alpha x + \beta y=0$ in $\C^2$ which does not lie  in the tangent cone of $(\Delta_p,0)$. 

 One may in turn define an element $h=a f + b g \in J$ with $a,b \in \OSO$ to be {\em general} iff $a(0) f + b(0) g$ defines a general element of the pencil $(C_{\alpha,\beta})$.





In fact, we define here, for {\em any} $m$-primary ideal $I$ in $\OSO$, a notion of  {\em general element} which has the following property~: take any pair $(f,g)$ of elements of~$I$ such that the ideal $J=(f,g)$ is a {\em reduction} of $I$ (see \S~\ref{sec-Samuel}),  then  the general elements of $J$ (in the ``pencil'' sense)  will be general elements of $I$, and conversely any general element of $I$ will be obtained as an element of such a reduction.

However, this will not be our first definition of the {\em general elements of} $I$ since we rather define them purely by their behaviour  on the normalized blow-up of~$I$ (cf. def.~\ref{def-elem-gene}). 

In a previous paper, we  proved that these elements are characterized by their Milnor number (theorem~\ref{thm-de-la-note}). Here, we focus on the applications of this result~:

 In \S~\ref{sec-two-special-cases}, we show how it covers both the study of limit of hyperplanes tangent to a normal surface, and the study of linear systems of plane curves, proving on one side a complement to a theorem by J. Snoussi, and on the other side an algebraic $\mu$-constant theorem for linear systems of plane curves (also obtained by other means by E. Casas).

In \S~\ref{sec-discri}, we prove the relation between our definition of general elements of $I$ and the one for pencils as claimed above. As a corollary, for two pencils $(f,g)$ and $(f',g')$ defining a  reduction of $I$, the general elements of both pencils have the same Milnor number, and the discriminants of the corresponding projections have the same multiplicity, which is apparently new.

\section{Geometry of a theorem by Samuel}
\label{sec-Samuel}

In this section only, we consider a germ  $(X,0)$ of complex analytic space with  arbitrary dimension $d$.
We let $\Oc:=\Oc_{X,0}$ be the corresponding local analytic ring. In fact, the content of this section would work with any local noetherian ring with infinite residue field (see e.g. \cite{Lipman-mult} or \cite{These} Chap.~2.3).

We recall that  an element $f\in \Oc$ is said to be {\em integrally dependant}  on an ideal $I$ of $\Oc$ if it satisfies an equation~:
\DebEq
$$f^n+a_1 f^{n-1}+\cdots+a_n=0,$$
\FinEq
with the condition $a_i\in I^i$ for all $i=1,\dots,n$.\par

The theory of integral dependance on ideals was initiated by O. Zariski (see \cite{Zar-Sam} Appendix 4) and under the influence of H. Hironaka was developed in the seminar \cite{LJ-Te}  where several characterizations are given. In the hands of B.~Teissier, it became a cornerstone in  the theory of equisingularity (see e.g. \cite{LaRabida} Chap.~1).
More recently, the theory  was extended to modules under the impulse of T. Gaffney (see the survey \cite{Gaffney-Massey-Trends}).

Let us just mention that the set $\overline{I}$ of the elements of $\Oc$ integrally dependant on $I$  is an ideal, called the integral closure of $I$ in $\Oc$, and that the definition of integral closure finds a natural expression on the blow-up $X_I$ of the germ $(X,0)$ along $I$ (see \cite{LaRabida}).

For the sake of simplicity, we restrict here to the case of a {\em reduced} germ $(X,0)$ (cf. \cite{These} loc. cit. for the general case). Then one  may take the normalization~$\overline{X_I}$  of the blow-up $X_I$, and 
following \cite{LaRabida} (Chap.~1, (1.3.6) et seq.), one proves that the equality $\overline{I}=\overline{J}$ of integral closures of ideals in $\Oc$ is equivalent to the equality~:
\DebEq
\begin{equation}
\label{eq-IOSI-JOSI}
I.\OXIB=J.\OXIB,
\end{equation}
\FinEq
\noindent for the corresponding sheaves on the normalized blow-up $\XIB$.

We now take $I$ to be an $m$-primary ideal of $\Oc$ i.e. containing a power $m^s$ of the maximal ideal of $\Oc$.

Denoting by $\bIb : \XIB \rightarrow (X,0)$ the normalized blow-up,  we write $D_1,\dots,D_s$ for the irreducible components of the reduced exceptional divisor $|(\overline{b_I})^{-1}(0)|$, and $v_{D_i}$ for the valuation along $D_i$.

Then we define (cf. \cite{Cras} d\'ef-prop.~1) an element $f\in I$ to be {\em $v$-superficial} if, and only if,
\DebEq
\begin{equation}
\label{eq-definit-v-super}
v_{D_i}(f)=v_{D_i}(I):=\inf \{ v_{D_i}(g), g \in I \}\; \mbox{for all  }\; i=1,\dots,s.
\end{equation}
Denoting by $D_f:=\sum_{i=1}^s v_{D_i}(f) D_i$, the total transform $(f)^*:=(f\circ\bIb)$ on $\XIB$ may be written as a sum of divisors~:
\DebEq
$$(f)^*=(f)'+D_f,$$
\noindent with $(f)'$ the strict transform of $f$ on $\XIB$.

The first part of the following proposition is an avatar of a theorem by P.~Samuel. The second part is the geometric version announced in the title~:

\begin{prop}
\label{prop-Samuel}
{\rm i)} Let $\Oc$ be a local noetherian ring of dimension $d$ with infinite residue field $\Oc/m$.
Let $I$ be an $m$-primary ideal of $\Oc$.
There exist a $d$-tuple $(f_1,\dots,f_d)$ of elements of $I$ such that the ideal $(f_1,\dots,f_d)$ is {\em a reduction} of $I$, i.e. has the same integral closure as $I$.

{\rm ii)} In our setting, let $\Oc$ be  the local ring of a reduced analytic germ $(X,0)$.  The $d$-tuple in i) are characterized by the two conditions that first, all the $f_i$ are $v$-superficial in $I$ and secondly, the intersection of their strict transforms  $(f_i)'$ on the exceptional divisor $\Dc$ verifies~:
$$(f_1)'\cap (f_2)'\cap \cdots \cap (f_d)'\cap \Dc=\emptyset.$$
We call such a $d$-tuple a {\em good} $d$-tuple of $v$-superficial elements in $I$.
\end{prop}

We will not give the proof here, but the reader should understand that ii) also easily gives the proof of  i) thanks to the characterization on (\ref{eq-IOSI-JOSI}) above. In fact, the same ``geometric proof'' works under the general  hypotheses of i) but one has to work on the non normalized blow-up (see~\cite{These} Chap.~2).

The original theorem by Samuel was formulated in terms of multiplicities (cf. \cite{Zar-Sam} Chap.~VIII thm.~22) so that it seems relevant to mention the following~:

\begin{prop}
\label{prop-v-super-e-transverse}
Let $\Oc$ be analytic local integral domain, and $I$ an $m$-primary ideal of $\Oc$.
The multiplicity $e(I/(f),\Oc/(f))=e(I,\Oc)$ if, and only if, $f$ is $v$-superficial.
\end{prop}

This result can  be deduced from a general formula for $e(I/(f),\Oc/(f))$ due to Flenner and Vogel in \cite{Fl-Vo} (for any noetherian local ring).

\section{General elements of an ideal}

From now on, we restrict ourselves to a two-dimensional {\em normal}  germ $(S,0)$.\footnote{For the elementary properties of normal surfaces we use here,  see \cite{Jawad} \S~2.6, and \cite{RES}.}
 
\begin{defi}
\label{def-elem-gene}
Let $\Oc$ be the local ring of a germ of normal surface $(S,0)$ and let  $I$ be an $m$-primary ideal of $\Oc$.
Adapting the notation from  section~\ref{sec-Samuel}, $\SIB$ denotes the normalized blow-up of $I$ on $(S,0)$.
We define an element  $f\in I$ to be {\em general} if, and only if, 
\begin{itemize}
\item[{\rm (i)}] $f$ is $v$-superficial in $I$ (cf. {\rm \S~\ref{sec-Samuel}~(\ref{eq-definit-v-super})}), 
\item[{\rm (ii)}] the strict transform $(f)'$ is a smooth curve transversal to the exceptional divisor $\Dc$ in $\SIB$,  which means that $(f)'$ does not  go through singular points either of  $\SIB$ or of $\Dc$ and that the intersection is transverse.
\end{itemize}
\end{defi}

Consider any resolution $r \: : \: X \rightarrow \SIB$ of the singularities of $\SIB$, {\em good} in the sense that, denoting  by $\pi=\bIb \circ r : X \rightarrow (S,0)$, the exceptional divisor $Z=\pi^{-1}(0)$  has only normal-crossing singularities.

Denote by $(f\circ \pi)=(f)'+Z_f$ the decomposition of the total transform of $(f)$ on $X$ into an exceptional (compact) part $Z_f$ 
and its strict transform denoted again $(f)'$.

Denoting $Z_I$ the divisor defined by $I.\Oc_X$ on $X$, we easily get the following~:

\begin{prop}
\label{prop-sur-la-resolution}
With the notation as above, $f\in I$ is general if, and only if, its total transform on $X$ is such  that~: 
\begin{itemize}
\item[$\alpha$)] its exceptional part is the generic one for the elements of $I$ i.e. $Z_f=Z_I$,
\item[$\beta$)] its strict transform is a (multi-germ of) smooth curves transversal to $Z$.
\end{itemize}
\end{prop}

As a corollary of this proposition, it is possible (either by a computation of Euler-Poincar\'e characteristic of covering spaces as indicated in \cite{Cras} \S~4, which followed \cite{GS},  or by an algebraic derivation from a Riemann-Roch formula as in \cite{Morales} 2.1.4)  to compute the Milnor number (in the sense of \cite{Bu-Gr}) of the complex curve defined by any general element $f\in I$.  We then  get~:
\begin{equation}
\label{eq-mu-general}
\mu(f)=\mu_I:=1-(Z_I.(Z_I-|Z_I|-K)),
\end{equation}
on any good resolution as defined before the proposition, where $|Z_I|$ (resp. $K$) denote the reduced divisor associated to $Z_I$ (resp. the numerically canonical cycle) and $(\: \cdot \:)$ denotes the intersection product  (see \cite{Cras}).

The main theorem in \cite{Cras} is the converse implication~:
\begin{thm}
\label{thm-de-la-note}
Let $(S,0)$ a germ of normal surface singularity, and $I$ an $m$-primary ideal of $\Oc_{S,0}$.
An element $f\in I$ is general in the sense of \ref{def-elem-gene} if, and only if, the Milnor number $\mu(f)$ has the value $\mu_I$ prescribed by formula (\ref{eq-mu-general}), which is also the minimum Milnor number for the elements of $I$.
\end{thm}

\begin{rem}
\label{rem-grace-a-Morales}
Thanks to the algebraic computation of the Milnor number for general elements which follows from \cite{Morales} (see before formula~(\ref{eq-mu-general})), theorem~\ref{thm-de-la-note} is proved without any topological argument, so that the proof fits to the setting of algebraic geometry over any algebraically  closed field of characteristic zero. 
\end{rem}

\section{Two special cases}

\label{sec-two-special-cases}
\subsection{The case when $(S,0)$ is arbitrary but $I=m$}

For a  germ $(S,0)$ of normal surface, given an embedding $(S,0)\subset (\C^N,0)$ defined by $N$ generators of the maximal ideal $m$ of $\Oc_{S,0}$,  we may consider the elements $f\in m$ as hypersurface sections of $S$.  From this point of view, J.~Snoussi studies in \cite{Jawad} what he called the {\em general hyperplanes} with respect to $(S,0)$. An hyperplane $H\ni 0$ in $\C^N$ is said to be {\em general} for $(S,0)$ if, and only if, it is not the limit of hyperplane tangents to the non singular locus of a small representant of $(S,0)$ in $\C^N$ (loc. cit. d\'ef. 2.2).  He then proves (loc. cit. thm.~4.2.)~:

\begin{thm}[Snoussi]
\label{thm-jawad}
If $(S,0)$ is a normal surface singularity embedded in $(\C^N,0)$, and if $H$ is an hyperplane {\em which does not contain an irreducible component of the tangent cone $C_{S,0}$ of $(S,0)$}, then $H$ is general if, and only if, the Milnor number $\mu(H\cap S,0)$ is minimum among the Milnor numbers of hyperplane sections of $(S,0)$.
 \end{thm}

From the definition of $v$-superficial elements given in \S~\ref{sec-Samuel}~(\ref{eq-definit-v-super}), it is clear that {\em the equation of an hyperplane $H$ defines a $v$-superficial element of the maximal ideal of $\Oc_{S,0}$  if, and only if, $H$ does not contain an irreducible component of $C_{S,0}$}. Hence, our theorem~\ref{thm-de-la-note} improves theorem~\ref{thm-jawad} as follows (see (ii))~:

\begin{cor}[of our theorem~\ref{thm-de-la-note}]
\label{cor-jawad-moi}
$ $

\noindent {\rm (i)} The equation of a general hyperplane in the sense of Snoussi  defines a general element of $m$ in the sense of definition~\ref{def-elem-gene}.  Conversely, if one takes a general element $f\in m$ and any embedding of $(S,0)\subset (\C^N,0)$ such that $f$ is induced by a coordinate function, $f=0$ defines a general hyperplane.

\noindent {\rm (ii)} In theorem~\ref{thm-jawad}, one may remove the hypothesis {\em ``$H$ does not contain an irreducible component of the tangent cone $C_{S,0}$ of $(S,0)$''} since theorem~\ref{thm-de-la-note} proved that elements with $\mu$ minimum {\em necessarily have this property}.\footnote{Note that this is exactly the tricky part of the argument in \cite{Cras}.}
\end{cor}

\subsection{The case when $(S,0)=(\C^2,0)$ and $I$ is arbitrary}

Since the definition of general element of an ideal given in def~\ref{def-elem-gene} was the same for an ideal $I$ and its integral closure $\bar{I}$ we consider only integrally closed ideals in the following discussion i.e. ideals such that $I=\bar{I}$.

These ideals were first studied by O. Zariski (see \cite{Zar-Sam} App. 5, where they are rather called {\em complete}) as the algebrization of Enriques'theory of clusters of points. For all this, we refer to the nice survey \cite{LJ}, and the book \cite{Casas}~: 
 
A cluster $K=(0_i,\nu_i)_i$ is a set of points $0_i$ {\em infinitely near $0$} i.e. lying above $0$ in a sequence of point blow-up starting from $(\C^2,0)$, with ascribed multiplicities~$\nu_i$. 
There is a one-to-one  correspondence between the integrally closed ideals of $(\C^2,0)$ and the clusters $(O_i,\nu_i)$ satisfying the so-called {\em proximity relations} of Enriques (see \cite{LJ} 5.1), also called {\em consistent clusters} in \cite{Casas} (p.~124).

For such a cluster $K=(0_i,\nu_i)_i$, the corresponding ideal $I_K$ is defined as the set of $f$ such that {\em virtual multiplicity} of the curve defined by $f$ at the point $0_i$ is at least $\nu_i$ (cf. loc. cit.).\footnote{In \cite{These} Chap.~1, we  explain how,  once $I_K$ is known, the virtual multiplicities of the elements of $I$ coincide with  the multiplicities of their weak transforms (cf. loc. cit.~1.1.6).}

Now $f\in I_K$ is said to {\em go sharply through $K$} if, and only if, $f$ goes through the $0_i$ with effective multiplicity equal to the $\nu_i$ {\em and} has no singular points outside~$K$ (cf. \cite{Casas} p.~127).

\begin{rem}
\label{rem-sharp-equising}
It is easy to see that two germs going sharply through $K$ are {\em equisingular} (cf. \cite{Casas} p.~127), in the sense of the well-known equisingularity theory of germs of plane curves.
\end{rem}

The careful study in \cite{LJ}, compared to our proposition~\ref{prop-sur-la-resolution}, yields~:

\begin{lem}
\label{lem-sharply-general}
For an integrally closed $m$-primary ideal $I$ of $\Oc_{\C^2,0}$, corresponding to a (consistent) cluster $K$, an element $f\in I$ is {\em general} in the sense of our def.~\ref{def-elem-gene} if, and only if, $f$ goes sharply through $K$ in the sense above.
\end{lem}
\begin{proof}
The proof of (ii)$\Leftrightarrow$(ii') in \cite{LJ} p.~360-361, gives exactly the equivalence between the fact that $f$ goes through the $O_i$ with effective multiplicity $\nu_i$ and the fact that $Z_f=Z_I$ on the minimal resolution of the blow-up of $(\C^2,0)$ along~$I$ (notation of prop.~\ref{prop-sur-la-resolution}).

Now the fact that $f$ has no singular points outside $K$ gives that the strict transform of $f$ on $S$ is transversal to the exceptional divisor by the argument of \cite{LJ}, proof of 6.1. (applied to each branch of $f$ corresponding to a simple ideal in the decomposition of $I$) . The converse is clear.
\end{proof}

With this, we get from our theorem~\ref{thm-de-la-note} and rem.~\ref{rem-sharp-equising} the following~:
\begin{cor}
\label{cor-mu-mini-donne-equising-C2}
For an integrally closed $m$-primary ideal $I$ of $\Oc_{\C^2,0}$, all the elements $f\in I$ such that 
$\mu(f)$ has the generic value $\mu_I$ are equisingular.
\end{cor}

Note that this  $\mu$-constant result for linear systems of plane curves is obtained without using topology (cf. rem.~\ref{rem-grace-a-Morales}). Another algebro-geometric proof of the same result is derived from the theory of clusters in \cite{Casas} \S~7.3. No such algebraic proof exists for the much more general theorem of L\^e (in \cite{Le-mu-constant}) on arbitrary (non-linear) family of germs of plane curves (cf. the remark on p.~361 in \cite{LaRabida}).

\section{General elements and discriminants}
\label{sec-discri}
We go back to our general setting i.e. $(S,0)$ is any germ of normal surface singularity, $I$ any  $m$-primary ideal of $\OSO$ and we take $(f,g)$ a good couple of $v$-superficial elements in $I$ (cf.  prop.~\ref{prop-Samuel}) so that $J=(f,g)$ is a reduction of $I$.
 
Let $p : (S,0) \rightarrow (\C^2,0)$ be the projection corresponding to $f,g$ as in the introduction, whose degree $\dgre(p)$ is by definition the multiplicity $e(f,g)=e(I)$.

Following Teissier (cf. \cite{Hunting}), one defines the critical space $(C_p,0)$ of $p$ by the ideal $I_{C_p}=F_0(\Omega_p)$ in $\Oc_{(S,0)}$,  where~$\Omega_p$ denotes the module of relative differentials,  and $F_0$ the zeroth Fitting ideal.
 Then,  denoting $\Oc_{C_p,0}=\OSO/I_{C_p}$, one constructs the discriminant space $(\Delta_p,0)$ as the image of $(C_p,0)$ by $p$,  defined in $(\C^2,0)$ by the ideal~:
\DebEq
$$I_{\Delta_p}:=F_0(p_*\Oc_{C_p,0}).$$
\FinEq
Now the space $(\Delta_p,0)$ may be both non-reduced at a generic point of one of its components, and have an embedded component at $0$. We denote $\Delta_{div}$ the divisorial part of $(\Delta_p,0)$ i.e. we do not consider the possible embedded component at $0$ (the reader will find more detail on all this in  \cite{RES} \S~3).
 
The following lemma was called L\^e-Greuel formula in \cite{RES} 3.9~:
\begin{lem}
\label{lem-Le-Greuel}
With the notation as above, 
for any  line $L\: :\: \alpha x +\beta y=0$             in  $\C^2$, denoting by   $(\;\cdot\;)$ the intersection number at $0$, we have the following equality~:\DebEq
\begin{equation}
\label{eq-Le-Greuel}
(\Delta_{div}\cdot L)_0=\mu(p^{-1}(L),0)+\dgre(p)-1,
\end{equation}
\FinEq
\noindent where $\mu$ is the Milnor number in the sense of {\rm  \cite{Bu-Gr}}.
\end{lem}

\begin{rem}
\label{rem-atteint-dans-pinceau}
From the definition~\ref{def-elem-gene} of  general elements applied to $J=(f,g)$, and Bertini's theorem, it is easy to see that for generic values of the numbers $(\alpha,\beta)\in \C^2$ the element $\alpha f+\beta g$ of the linear pencil defined by $f,g$  {\em is a general element of the ideal $J$} (cf. also \cite{Cras}).
\end{rem}

Now with the formula~(\ref{eq-Le-Greuel}) above, one deduces the following~:

\begin{cor}[of theorem~\ref{thm-de-la-note}]
\label{cor-carac-element-general}
Let $J=(f,g)$ be an $m$-primary ideal of $\Oc_{S,0}$. The elements $\alpha f+\beta g$ with $(\alpha,\beta)\in \C^2$ which are {\em general elements of $J$} (in the sense of def.~\ref{def-elem-gene}) are exactly the inverse-images $p^{-1}(L)$ of the lines $L~: \alpha x+\beta y=0$ transversal to the discriminant in lemma~\ref{lem-Le-Greuel}.
\end{cor}
\begin{proof}
The formula~(\ref{eq-Le-Greuel}) gives the equivalence between minimal Milnor number in the pencil and minimal intersection number $(\Delta_{div}\cdot L)$; by remark~\ref{rem-atteint-dans-pinceau}  we already know that the minimum Milnor number in $J$ is obtained by elements of the pencil, and we conclude by theorem~\ref{thm-de-la-note}.
\end{proof}

Now, we may also compare two projections $p=(f,g)$ and $p'=(f',g')$ such that the corresponding ideals $(f,g)$ and $(f',g')$ have the same integral closure $I$. We then know that $\dgre(p)=\dgre(p')=e(I)$ and from the foregoing, the generic Milnor numbers in the two pencils defined by $(f,g)$ and $(f',g')$ are  both equal to the same number $\mu_I$ as defined in formula (\ref{eq-mu-general}). Then  formula~(\ref{eq-Le-Greuel}) yields~:

\begin{cor}
\label{cor-mult-discri}
Let $(f,g)$ and $(f',g')$ be two $m$-primary ideals on $\Oc_{S,0}$ having the same integral closure $I$, 
the multiplicity of the discriminant of the projections $p$ defined by $f,g$ and $p'$ defined by $f',g'$ onto $\C^2$ are the same, equal to~:
\DebEq
\begin{equation}
\label{eq-Delta-I}
e(\Delta_p,0)=\mu_I+e(I)-1.
\end{equation}
\end{cor}

In the special case of any  projection $p=(f,g)$ with $\dgre(p)$ equal to the multiplicity $e(S,0)$ of    the germ $(S,0)$ (which is by definition the multiplicity $e(m)$), by a theorem of Rees (cf. \cite{LaRabida} p. 340), $(f,g)$ is a reduction of $m$.  Hence, for any such projection we have the following formula~:
\DebEq
\begin{equation}
\label{equation-Delta-m}
e(\Delta_p,0)=\mu^2+\mu^1,
\end{equation}
\FinEq
where $\mu^2$ is the generic Milnor number $\mu_m$ in $m$, and $\mu^1=e(S,0)-1$. The notation $\mu^i$ follows Teissier (cf. \cite{LaRabida} ex.~2.2 p.~423)  where the $\mu^i(X,0)$ denote in general the Milnor number of the intersection of a hypersurface $(X,0)$  with a ``general enough'' linear subspace  of dimension $i$.  Here our $(S,0)$ is no longer an hypersurface so that the Milnor number is in the sense of \cite{Bu-Gr}.  Note also that even in the case of an hypersurface of $\C^3$,  we obtain our formula (\ref{equation-Delta-m}) for a broader class of projections than the corresponding formula in \cite{LaRabida} loc. cit.

\medskip

{\footnotesize
\noindent{\bf Acknowledgement}
The idea of making more explicit the consequences of our previous work \cite{Cras} partly arose from the conversations I had at the CIRM with several ``Franco-Japanese'' mathematicians. I would like to thank them for their interest and the organizers  for such a nice opportunity.}

\bigskip
\noindent
C.M.I., 39 rue F. Joliot-Curie,\\
13453 Marseille Cedex 13\\
France\\
email~: bondil@cmi.univ-mrs.fr

\end{document}